 \documentclass[final,12pt,3p,times]{elsarticle}

\usepackage{amssymb}
\usepackage{amsthm}
\usepackage{hyperref}
\usepackage{amsmath,amssymb,amsopn,amsfonts,mathrsfs,amsbsy,amscd}
\usepackage{longtable}
\usepackage{mathtools}
\journal{~}

\newtheorem{theorem}{Theorem}
\newtheorem{definition}[theorem]{Definition}
\newtheorem{proposition}[theorem]{Proposition}

\newtheorem{corollary}[theorem]{Corollary}

\newtheorem{example}[theorem]{Example}

\newtheorem{remark}[theorem]{Remark}
\numberwithin{equation}{section}

\newtheorem{oquestion}{Question}

\newcommand{\dem}{\noindent{\bf Proof. }}
\DeclareMathAlphabet{\mathpzc}{OT1}{pzc}{m}{it}
\newcommand*{\mydot}{\mathrel{\scalebox{0.70}{$\odot$}}}
\begin{document}

\begin{center}
{\Large  Symmetric Zinbiel superalgebras}\footnote{
The work is supported by FCT   UIDB/00212/2020 and UIDP/00212/2020; RFBR 20-01-00030.
The authors thank Pasha Zusmanovich for some useful discussions.} 
\end{center} 

\begin{frontmatter}
\author[label1]{Sa{\"i}d Benayadi}
\address[label1]{Universit\'e de Lorraine\\Laboratoire de Math\'ematiques  IECL
UMR CNRS 7502\\3 rue Augustin Fresnel, BP 45112,  F-57073 Metz Cedex 03, France\\
{e-mail: said.benayadi@univ-lorraine.fr}}

\author[label2]{Ivan Kaygorodov}
\address[label2]{Centro de Matemática e Aplicações,  Universidade da Beira Interior, Covilhã, Portugal\\

Moscow Center for Fundamental and Applied Mathematics, Moscow,   Russia\\
 
{e-mail: kaygorodov.ivan@gmail.com}}

\author[label3]{Fahmi Mhamdi}
\address[label3]{Facult\'e des Sciences de Gafsa\\
D\'epartement de Math\'ematiques\\
Zarroug 2112, Tunisie\\
{e-mail: fahmi.mhamdi@gmail.com}}

\end{frontmatter}

 \
 
 \

\noindent {\bf Abstract.}
{\it The notion of symmetric Zinbiel superalgebras is introduced.
We prove that the nilpotency index of a symmetric Zinbiel superalgebra is not greater than 4 and describe two-generated symmetric Zinbiel algebras and odd generated superalgebras.
We discuss identities of mono and binary symmetric Zinbiel and Leibniz algebras. 
It is proven that each quadratic Zinbiel algebra is $2$-step nilpotent.
Also, we study double extensions of symmetric Zinbiel algebras.

}

\ 


\ 


 \section*{Introduction}

Loday introduced a class of symmetric operads generated by one bilinear operation subject to one relation making each left-normed product of three elements equal to a linear combination of right-normed products:
   \begin{center} $(a_1a_2)a_3=\sum\limits_{\sigma\in \mathbb{S}_3} x_{\sigma} a_{\sigma(1)}(a_{\sigma(2)}a_{\sigma(3)});$\end{center}
such an operad is called a parametrized one-relation operad. For a particular choice of parameters $\{x_{\sigma}\}$, this operad is said to be regular if each of its components is the regular representation of the symmetric group; equivalently, the corresponding free algebra on a vector space $V$ is, as a graded vector space, isomorphic to the tensor algebra of $V$. 
Bremner and Dotsenko classified, over an algebraically closed field of characteristic zero, all regular parametrized one-relation operads. In fact, they proved that each such operad is isomorphic to one of the following five operads: 
the left-nilpotent operad, the associative operad, the Leibniz operad, the Zinbiel operad, and the Poisson operad \cite{bredo}. 
An algebra $\mathcal{A}$ is called a (left) {\it Zinbiel algebra} if it satisfies the identity 
\begin{center}$(xy)z=x(yz+zy).$\end{center}
	Zinbiel algebras were introduced by Loday in \cite{loday}.
Under the Koszul duality, the operad of Zinbiel algebras is dual to the operad of Leibniz algebras. 
Zinbiel algebras are also known as pre-commutative algebras \cite{ pasha}
and chronological algebras \cite{Kawski}.
A Zinbiel algebra is equivalent to a commutative dendriform algebra \cite{comdend}.
It plays an important role in the definition of pre-Gerstenhaber algebras.
The variety of Zinbiel algebras is a proper subvariety in the variety of right commutative algebras.
Each  Zinbiel algebra with the commutator multiplication gives a Tortkara algebra \cite{dzhuma}, which has sprung up in unexpected areas of mathematics \cite{tortnew1,tortnew2}.
Recently, the notion of matching Zinbiel algebras was introduced in \cite{matchzin}.
Zinbiel algebras also  appeared in a study of rack cohomology \cite{rack}, 
number theory \cite{Chapoton21} and 
in  a  construction of a Cartesian differential category  \cite{ip21}.
In recent years, there has been a strong interest in the study of Zinbiel algebras in the algebraic and the operad context  
\cite{  centr3zinb,    cam13, dok, ims21,    dzhuma,    dzhuma5, dzhuma19, matchzin,  kppv, 34,mukh,  anau, ualbay}. 

Free Zinbiel algebras were shown to
be precisely the shuffle product algebra \cite{34}, 
which is under a certain interest until now \cite{cl}.
Naurazbekova proved that, over a field of characteristic zero, free Zinbiel algebras are the free associative-commutative algebras (without unity) with respect to the symmetrization multiplication and their free generators are found; also she constructed examples of subalgebras of the two-generated free Zinbiel algebra that are free Zinbiel algebras of countable rank \cite{anau}.
Nilpotent algebras play an important role in the class of Zinbiel algebras.
So, 
Dzhumadildaev and  Tulenbaev proved that each complex finite-dimensional Zinbiel algebra is nilpotent  \cite{dzhuma5}
and recently  Towers proved that each finite-dimensional Zinbiel algebra is nilpotent \cite{zinb5};
Naurazbekova 
and Umirbaev proved that in characteristic zero any proper subvariety of the variety of Zinbiel algebras is nilpotent \cite{ualbay}.
Finite-dimensional Zinbiel algebras with a ``big'' nilpotency index are classified in \cite{adashev,cam13}.
Central extensions of three-dimensional Zinbiel algebras were calculated in \cite{centr3zinb} and of filiform Zinbiel algebras in \cite{cam20}.
The full system of degenerations of complex four-dimensional Zinbiel algebras is given in \cite{kppv}.

The present paper is about symmetric (left and right) Zinbiel (super)algebras.
This variety of algebras is dual to symmetric Leibniz algebras 
(about symmetric Leibniz algebras see in \cite{BMO} and references therein) and it is very thin.
As we have proven each symmetric Zinbiel algebra is a  $2$-step nilpotent or $3$-step nilpotent algebra.
We also provide a classification of two-generated symmetric Zinbiel algebras (Theorem \ref{thalg2gen}) and symmetric Zinbiel superalgebras with two odd generators (Proposition \ref{twood}).
The next part of the paper is devoted to a study of mono and binary symmetric Zinbiel algebras. Next, we give the description of quadratic Zinbiel algebras.
It has been proven that each quadratic Zinbiel algebra is $2$-step nilpotent (Proposition \ref{quadra}).
The notion of double extensions of quadratic Zinbiel algebras is introduced and described in the last section.

\section{Definitions and preliminary results}\label{sect1}

Let $V=V_{\bar{0}}\oplus V_{\bar{1}}$ be a $\mathbb{Z}_{2}$-graded vector space over the field $\mathbb{K}$. An element $x$ in $V$ is called homogeneous if $x\in V_{\bar{0}}$ or $x\in V_{\bar{1}}$. Throughout this paper, all elements are supposed to be homogeneous unless otherwise stated. For a homogeneous element $x$ we shall use the standard notation $\vert x\vert\in \mathbb{Z}_{2}=\{\bar{0},\bar{1}\}$ to indicate its degree, i.e. whether it is contained in the even part ($\vert x\vert=\bar{0}$) or in the odd part ($\vert x\vert=\bar{1}$). All  superalgebras considered in this paper are finite-dimensional.

Let $V$ be a $\mathbb{K}-$vector space and $\gamma:= {\wedge V}$ a Grassmann (or exterior) algebra of $V$. We know that $\Gamma:= \oplus_{i\in \mathbb{Z}}\wedge^iV= \Gamma_{\bar{0}}\oplus \Gamma_{\bar{1}}$ is a
$\mathbb{Z}_{2}-$graded associative algebra, where $ \Gamma_{\bar{0}}:= \oplus_{i\in \mathbb{Z}}\wedge^{2i}V$ and $ \Gamma_{\bar{1}}:= \oplus_{i\in \mathbb{Z}}\wedge^{2i+1}V$, such that
$X_\alpha X_\beta= (-1)^{\alpha\beta }X_\beta X_\alpha, \,\, \forall \, (X_\alpha,X_\beta)\in \Gamma_\alpha\times \Gamma_\beta.$

Let $(\mathcal{A}=\mathcal{A}_{\bar{0}}\oplus \mathcal{A}_{\bar{1}}, \cdot)$ be a superalgebra and $\Gamma(\mathcal{A})$ its Grassmann enveloping algebra which is a subalgebra of $\mathcal{A}\otimes \Gamma$ given by $\Gamma(\mathcal{A})=\mathcal{A}_{\bar{0}}\otimes \Gamma_{\bar{0}}\oplus \mathcal{A}_{\bar{1}}\otimes \Gamma_{\bar{1}}$. Let us assume now that $\mathfrak{B}$ is a homogeneous variety of algebras. Then, $\mathcal{A}$ is said to be a $\mathfrak{B}$-superalgebra if $\Gamma(\mathcal{A})$ belongs to $\mathfrak{B}$. So, $\mathcal{A}=\mathcal{A}_{\bar{0}}\oplus \mathcal{A}_{\bar{1}}$ is a left (resp. right) Zinbiel superalgebra if $\Gamma(\mathcal{A})$ is a left (resp. right) Zinbiel algebra. Consequently, we get the following definition.
\begin{definition}\label{d11}
Let $\mathcal{A}=\mathcal{A}_{\bar{0}}\oplus \mathcal{A}_{\bar{1}}$ be a $\mathbb{Z}_{2}$-graded vector space and let $``\cdot": \mathcal{A} \times \mathcal{A} \rightarrow \mathcal{A}$ be a bilinear map on $\mathcal{A}$ such that $\mathcal{A}_{i} \cdot \mathcal{A}_{j}\subset\mathcal{A}_{i+j},~\forall i,j\in \mathbb{Z}_{2}$.
\begin{enumerate}
\item[1.]
$\mathcal{A}$ is called a left (resp. right) Zinbiel superalgebra if, for all $x,y,z\in \mathcal{A}_{ \bar{0}} \cup \mathcal{A}_{\bar{1} },$
\begin{equation}\label{IdentLeftZinb}
(xy)z=x(yz)+(-1)^{\vert y\vert \vert z\vert} x(zy),
\end{equation}
\vspace*{-0.35cm}
\begin{equation}\label{IdentRightZinb}
({\rm resp.}~~
x(yz)=(xy)z+(-1)^{\vert x\vert \vert y\vert}(yx)z),
\end{equation}
or, equivalently,
$(x,y,z)=(-1)^{\vert y\vert \vert z\vert}x(zy)$
( \ resp. $(x,y,z)=-(-1)^{\vert x\vert \vert y\vert}(yx)z \ )$,
where $(x,y,z)=(xy)z-x(yz)$ is the associator associated to $``\cdot"$.
\item[2.] $\mathcal{A}$ is called a symmetric Zinbiel superalgebra if it is both a left and a right Zinbiel superalgebra.
\end{enumerate}
\end{definition}

As usual, for $x  \in \mathcal{A}_{ \bar{0}} \cup \mathcal{A}_{\bar{1} },$ we define the corresponding endomorphism of $\mathcal{A}$ by $L_{x}(y)=xy$  (resp. $R_{x}(y)=(-1)^{\vert x\vert \vert y\vert}yx$), $\forall y\in \mathcal{A}_{ \bar{0}} \cup \mathcal{A}_{\bar{1} }$, which is called the left (resp. right) multiplication by $x$.

The following proposition can be verified by direct calculation.

\begin{proposition}\label{CondSym}
Let $(\mathcal{A}, \cdot)$ be a symmetric Zinbiel superalgebra.
Then,
\begin{enumerate}
    \item $(\mathcal{A},\cdot)$ is  a LR-superalgebra (LR-superalgebras are also known as bicommutative superalgebras), i.e. it satisfies
    $$(xy)z=(-1)^{\vert y\vert \vert z\vert}(xz)y \mbox{ and }x(yz)=(-1)^{\vert x\vert \vert y\vert}y(xz), \mbox{ for all }x,y,z\in \mathcal{A}_{ \bar{0}} \cup \mathcal{A}_{\bar{1} },$$

\item $(\mathcal{A}, \cdot)$ is an anti-flexible superalgebra, i.e. it satisfies
$$(x,y,z)=(-1)^{\vert x\vert \vert y\vert+\vert z\vert (\vert x\vert + \vert y\vert)}(z,y,x), \mbox{ for all }x,y,z\in \mathcal{A}_{ \bar{0}} \cup \mathcal{A}_{\bar{1} },$$

\item $(\mathcal{A}, \cdot)$ satisfies
$$(xy)z=-(-1)^{\vert x\vert (\vert y\vert + \vert z\vert)}y(zx), \mbox{ for all }x,y,z\in \mathcal{A}_{ \bar{0}} \cup \mathcal{A}_{\bar{1} },$$

\item $(\mathcal{A}, \cdot)$ satisfies
$$(xy)z=-(-1)^{\vert x\vert\vert y\vert +\vert y\vert\vert z\vert + \vert z\vert\vert x\vert}z(yx), \mbox{ for all }x,y,z\in \mathcal{A}_{ \bar{0}} \cup \mathcal{A}_{\bar{1} }.$$

\end{enumerate}
\end{proposition}

In the following, we recall some results about representations of left and right Zinbiel superalgebras.
\begin{definition} \label{d16}
Let $(\mathcal{A}, \cdot)$ be a non-associative superalgebra, $\mathcal{V}$ be a $\mathbb{Z}_{2}$-graded vector space and $r,l:\mathcal{A}\rightarrow End(\mathcal{V})$ be two even linear maps. If $\mathcal{A}$ is a left (resp. right) Zinbiel superalgebra, then we say that $(r,l)$ is a left (resp. right) representation of $\mathcal{A}$ in $\mathcal{V}$ if for all $x,y \in \mathcal{A}_{\vert x\vert}\times \mathcal{A}_{\vert y\vert}$:
$$l(xy)=l(x)l(y)+l(x)r(y);~~r(x)r(y)=r(xy)+(-1)^{\vert x\vert \vert y\vert}r(yx);~~l(x)l(y)=
(-1)^{\vert x\vert \vert y\vert}[r(y),l(x)],$$
$$\big({\rm resp.}~~l(x)l(y)=l(xy)+(-1)^{\vert x\vert \vert y\vert}l(yx);~~~~l(x)r(y)=r(xy)=(-1)^{\vert x\vert \vert y\vert}(r(y)r(x)+r(y)l(x))\big).$$
\end{definition}
\begin{example} \label{adjrep}\rm
Let $(\mathcal{A}, \cdot)$ be a left (resp. right) Zinbiel superalgebra. Consider the even maps $L:\mathcal{A}\rightarrow End(\mathcal{A})$ and $R:\mathcal{A}\rightarrow End(\mathcal{A})$ defined by: $L(x):=L_{x},~R(x):=R_{x},\forall x\in \mathcal{A}$. Then, $(R,L)$ is a left (resp. right) representation of $\mathcal{A}$ in $\mathcal{A}$ called the left (resp. right) adjoint representation of $\mathcal{A}$.
\end{example}
\begin{proposition} \label{p18}
Let $(\mathcal{A},\cdot)$ be a left (resp. right) Zinbiel superalgebra and $r,l:\mathcal{A}\rightarrow End(\mathcal{V})$ be two even linear maps. Then, the $\mathbb{Z}_{2}$-graded vector space $\overline{\mathcal{A}}=\mathcal{A}\oplus \mathcal{V}$ endowed with the product defined by: $(x+u)  (y+v)=x  y+l(x)(v)+(-1)^{\vert x\vert \vert y\vert}r(y)(u),~\forall (x+u,y+v)\in \overline{\mathcal{A}}_{\vert x\vert}\times \overline{\mathcal{A}}_{\vert y\vert}$, is a left (resp. right) Zinbiel superalgebra if and only if $(r,l)$ is a left (resp. right ) representation of $\mathcal{A}$ in $\mathcal{V}$.
\end{proposition}
\begin{proposition} \label{corep}
Let $(\mathcal{A}, \cdot)$ be a left (resp. right) Zinbiel superalgebra and $(R,L)$ be the left (resp. right) adjoint representation of $\mathcal{A}$. Let us consider the even linear maps $L^*,R^*:\mathcal{A}\rightarrow End(\mathcal{A}^*)$ defined by:\[L^*(x)(f)=(-1)^{\vert f\vert \vert x\vert}f\circ R(x),~~~~R^*(x)(f)=(-1)^{\vert f\vert \vert x\vert}f\circ L(x),~~\forall x\in \mathcal{A}_{\vert x\vert},f\in \mathcal{A}^*_{\vert f\vert}~.\]$(R^*,L^*)$ is a left (resp. right) representation of $\mathcal{A}$ in $\mathcal{A}^*$ if and only if $\mathcal{A}$ is $2$-step nilpotent.
\end{proposition}
\dem Suppose that $(\mathcal{A},\cdot)$ is a left Zinbiel superalgebra. Then, for all
$x\in \mathcal{A}_{\vert x\vert},y\in \mathcal{A}_{\vert y\vert},z\in \mathcal{A}_{\vert z\vert}$ and $f\in \mathcal{A}^*_{\vert f\vert}$, a simple calculation proves that $L^*(xy)(f)(z)=\big(L^*(x)L^*(y)+L^*(x)R^*(y)\big)(f)(z)$
if and only if
\begin{equation}\label{IDE1}
y(zx)=-(-1)^{\vert y\vert \vert z\vert}z(yx).
\end{equation}
Moreover,
$R^*(x)R^*(y)(f)(z)=\big(R^*(xy)+(-1)^{\vert x\vert \vert y\vert}R^*(yx)\big)(f)(z)$ if and only if
\begin{equation}\label{IDE2}
y(zx)=-(-1)^{\vert x\vert (\vert y\vert+\vert z\vert)}(xy)z.
\end{equation}
Furthermore, according to \eqref{IDE1}, $L^*(x)L^*(y)(f)(z)=\big((-1)^{\vert x\vert \vert y\vert}R^*(y)L^*(x)-L^*(x)R^*(y)\big)(f)(z)$ if and only if
\begin{equation}\label{IDE3}
z(xy)=(-1)^{\vert y\vert (\vert x\vert+\vert z\vert)}(yz)x.
\end{equation}
By \eqref{IDE2} and \eqref{IDE3}, we conclude that $(R^*,L^*)$ is a left representation of $\mathcal{A}$ in $\mathcal{A}^*$ if and only if $\mathcal{A}$ is $2$-step nilpotent. By the same way, we prove the
result for right Zinbiel superalgebras.\hfill $\square$\\

\quad In the following, we introduce the representation of symmetric Zinbiel superalgebras and give some related results.
\begin{definition} \label{d112}
Let $(\mathcal{A}, \cdot)$ be a symmetric Zinbiel superalgebra, $\mathcal{V}$ be a $\mathbb{Z}_{2}$-graded vector space and  $r,l:\mathcal{A}\rightarrow End(\mathcal{V})$ be two even linear maps. Then, we say that $(r,l)$ is a representation of $\mathcal{A}$ in $\mathcal{V}$ if $(r,l)$ is a left and a right representation of $\mathcal{A}$ in $\mathcal{V}$. We denote by $Rep(\mathcal{A},\mathcal{V})$ the set of all representations of $\mathcal{A}$ in $\mathcal{V}$.
\end{definition}
\begin{example} \label{e113}\rm
If $(\mathcal{A},\cdot)$ is a symmetric Zinbiel superalgebra, then:
\begin{enumerate}
\item[1.] $(R,L)$ (see Example \ref{adjrep}) is a representation of $\mathcal{A}$ in $\mathcal{A}$ called the adjoint representation of $\mathcal{A}$.
\item[2.] $(R^*,L^*)$ (see Proposition \ref{corep}) is a representation of $\mathcal{A}$ in $\mathcal{A}^*$ called the co-adjoint representation of $\mathcal{A}$.
\end{enumerate}
\end{example}
\begin{remark} \label{SLSR}
Let $(\mathcal{A},\cdot)$ be a symmetric Zinbiel superalgebra, $\mathcal{V}$ be a $\mathbb{Z}_{2}$-graded vector space and $r,l:\mathcal{A}\rightarrow End(\mathcal{V})$ be two even linear maps. Then, the $\mathbb{Z}_{2}$-graded vector space $\overline{\mathcal{A}}=\mathcal{A}\oplus \mathcal{V}$ endowed with the product defined by: \begin{center}$(x+u)(y+v)=xy+l(x)(v)+(-1)^{\vert x\vert \vert y\vert}r(y)(u),~\forall (x+u,y+v)\in \overline{\mathcal{A}}_{\vert x\vert}\times \overline{\mathcal{A}}_{\vert y\vert}$, 
\end{center} is a symmetric Zinbiel superalgebra if and only if $(r,l)\in Rep(\mathcal{A},\mathcal{V})$.
\end{remark}

\section{Characterizations of symmetric Zinbiel superalgebras}\label{sect2}
This section will be devoted to the characterizations of symmetric Zinbiel superalgebras.
In this section, we will consider algebras over the complex field.

\begin{definition}
Let $(\mathcal{A},\cdot)$ be a superalgebra, then
\begin{itemize}
    \item $(\mathcal{A},\cdot)$ is a $2$-step nilpotent superalgebra if for all  elements $x,y,z\in \mathcal{A}_{ \bar{0}} \cup \mathcal{A}_{\bar{1} }$ we have $(xy)z=x(yz)=0,$ but the product of the algebra is nonzero.

    \item $(\mathcal{A},\cdot)$ is a $3$-step nilpotent superalgebra if for all  elements $x,y,z,t\in \mathcal{A}_{ \bar{0}} \cup \mathcal{A}_{\bar{1} }$ we have \begin{center}$(xy)(zt)=x((yz)t)=x(y(zt))=((xy)z)t=(x(yz))t=0,$ 
    \end{center} but the algebra is non-$2$-step nilpotent.

\end{itemize}

\end{definition}

\begin{example}
Each $2$-step nilpotent superalgebra is a symmetric Zinbiel superalgebra.
\end{example}

\begin{theorem} \label{t510}
If $(\mathcal{A},\cdot)$ is a nonzero non-$2$-step nilpotent symmetric Zinbiel superalgebra, then
\begin{enumerate}
\item $(\mathcal{A},\cdot)$ is a $3$-step nilpotent superalgebra, i.e. it is a central extension of a $2$-step nilpotent superalgebra;
\item  $x^3=0$ for each $x \in \mathcal{A}_{\vert x\vert}.$
\end{enumerate}
\end{theorem}
\dem
For the   first part of the theorem we consider the following relations:
\begin{longtable}{lcl}
$(xy)(zt)$& $=$&
$(-1)^{\vert z\vert(\vert x \vert+\vert y\vert)} z((xy)t) =
-(-1)^{(\vert z\vert+\vert t\vert)(\vert x \vert+ \vert y \vert)+ \vert x \vert \vert y \vert } z(t(yx)),$\\

$(xy)(zt)$& $=$&
$ -(-1)^{(\vert x\vert+ \vert y \vert)(\vert z \vert+\vert t \vert)+ \vert z \vert \vert t\vert }(tz)(xy)=
-(-1)^{\vert t \vert(\vert x \vert+ \vert y \vert+ \vert z \vert)} (t(xy))z =
(-1)^{(\vert z\vert+\vert t\vert)(\vert x \vert+ \vert y \vert)} z(t(xy)).$
\end{longtable}
Hence,
\[0= z(t(xy +  (-1)^{  \vert x \vert \vert y \vert } yx))= z((tx)y),\]
which gives that all products of $4$ elements are equal to zero and
$(\mathcal{A},\cdot)$ is a $3$-step nilpotent superalgebra.

Now, for the second part of the theorem, we are using \eqref{IdentLeftZinb} and \eqref{IdentRightZinb} and substituting $y=z=x\in \mathcal{A}_{\vert x \vert}$. We get:
\begin{enumerate}[(a)]
\item If $x\in \mathcal{A}_{\bar{0}}$, then $x^{2}x=2xx^{2}$ and $xx^{2}=2x^{2}x$, which implies that $x^{2}x=0=xx^{2}$.
So, $x^{3}=x^{2}x=xx^{2}=0$.
\item If $x\in \mathcal{A}_{\bar{1}}$, then $x^{2}x=0=xx^{2}$. So, $x^{3}=0$.
\end{enumerate}
 \hfill $\square$

\begin{corollary}
If $(\mathcal{A}, \cdot)$ is a $d$-generated symmetric Zinbiel superalgebra, then 
$dim(\mathcal{A}, \cdot) \leq -d+d^2+2d^3+d^4.$
\end{corollary}

\dem
It is easy to see that the maximal dimension   of a central extension $\hat{\mathcal{A}}$ of the $d$-dimensional zero product algebra (with a basis $\{e_1, \ldots, e_d\}$)
is $d+d^2.$
Let say that the last algebra has a basis $\{e_1, \ldots, e_d, e_{1,1}, \ldots, e_{d,d} \},$ such that $e_{i,j}=e_ie_j.$
${\mathcal{A}}$ is a central extension of $\hat{\mathcal{A}}.$
Note now, that $x^3=0$ for each $x \in \mathcal{A},$
hence $e_ie_{i,i}=e_{i,i}e_i=0$ and ${\rm dim} \ \hat{\mathcal{A}}^2=d^2.$
It is easy to see that the dimension of the cohomology space of 
 $\hat{\mathcal{A}}$ is not greater than $(d+d^2)^2-d^2-2d.$
Hence, the dimension of a central extension of a central extension of the  $d$-generated zero product algebra 
is  not greater than $d+d^2+(d+d^2)^2-d^2- 2d=-d+d^2+2d^3+d^4.$
\hfill $\square$

\begin{corollary}
If $(\mathcal{A}, \cdot)$ is a one-generated symmetric Zinbiel superalgebra, then there are only
two following opportunities
\begin{enumerate}
    \item $dim(\mathcal{A}, \cdot)=(2,0)$ and the multiplication table is given by $e_1^2=e_2.$
    \item $dim(\mathcal{A}, \cdot)=(1,1),$ $\mathcal{A}_{\bar{0}}= \langle e_1 \rangle,$ $\mathcal{A}_{\bar{1}}= \langle e_2 \rangle,$ and the multiplication table is  given by $e_2^2=e_1.$

\end{enumerate}

\end{corollary}

\begin{theorem}\label{thalg2gen}
If $(\mathcal{A}, \cdot)$ is a two-generated symmetric Zinbiel algebra, then there are only
 following opportunities:

\begin{enumerate}
    \item
if    $(\mathcal{A}, \cdot)$ is a $2$-step nilpotent algebra, then $ 3 \leq dim(\mathcal{A}, \cdot) \leq 6$ and
\begin{enumerate}
    \item if $dim(\mathcal{A}, \cdot)=3,$ then $(\mathcal{A}, \cdot)$ is isomorphic to one of the following algebras
$$\begin{array}{ll llll}
\mathfrak{N}^3_1 &:& e_1 e_1 = e_2 &  e_1 e_2=0 & e_2 e_1=0  & e_2 e_2=0  \\
\mathfrak{N}^3_2 &:& e_1 e_1 = e_3 & e_1 e_2=0 & e_2 e_1=0  & e_2 e_2=e_3 \\
\mathfrak{N}^3_3 &:& e_1 e_1 = 0 & e_1 e_2=e_3 & e_2 e_1=-e_3  & e_2 e_2=0 \\
\mathfrak{N}^{3}_{4,\lambda } &:& e_1 e_1 = \lambda e_3 & e_1 e_2=0 & e_2 e_1=e_3  & e_2 e_2=e_3 \\
\end{array}$$
    \item if $dim(\mathcal{A}, \cdot)=4,$ then $(\mathcal{A}, \cdot)$ is isomorphic to one of the following algebras

$$\begin{array}{ll llll}
\mathfrak{N}^4_1  &:& e_1 e_1 = 0 & e_1 e_2=e_3 & e_2 e_1=e_4  & e_2 e_2=0 \\
\mathfrak{N}^4_2&:& e_1 e_1 = e_3 & e_1 e_2=0 & e_2 e_1=e_4  & e_2 e_2=0 \\
\mathfrak{N}^4_3 &:& e_1 e_1 = e_3 & e_1 e_2=0 & e_2 e_1=0  & e_2 e_2=e_4 \\
\mathfrak{N}^4_4 &:& e_1 e_1 = e_3 & e_1 e_2=e_3 & e_2 e_1=e_4  & e_2 e_2=e_3 \\
\mathfrak{N}^4_{5,\lambda} &:& e_1 e_1 = e_3 & e_1 e_2=0 & e_2 e_1=e_4  & e_2 e_2=e_3+\lambda e_4 \\
\mathfrak{N}^4_{6,\lambda } &:& e_1 e_1 = e_3 & e_1 e_2=e_4 & e_2 e_1=\lambda e_4  & e_2 e_2=0 \\
\end{array}$$

\item if $dim(\mathcal{A}, \cdot)=5,$ then $(\mathcal{A}, \cdot)$ is isomorphic to one of the following algebras

$$\begin{array}{ll llll}
\mathfrak{N}^5_1 &:& e_1 e_1 = 0 & e_1 e_2=e_3 & e_2 e_1=e_4  & e_2 e_2=e_5 \\
\mathfrak{N}^5_2 &:& e_1 e_1 = e_5 & e_1 e_2=e_3 & e_2 e_1=e_4  & e_2 e_2=e_5 \\
\mathfrak{N}^5_3 &:& e_1 e_1 = e_3 & e_1 e_2=e_4 & e_2 e_1=e_4  & e_2 e_2=e_5 \\
\mathfrak{N}^5_{4, \lambda} &:& e_1 e_1 = e_3+\lambda e_5 & e_1 e_2=e_3 & e_2 e_1=e_4  & e_2 e_2=e_5 \\
\end{array}$$

\item if $dim(\mathcal{A}, \cdot)=6,$ then $(\mathcal{A}, \cdot)$ is isomorphic to  the following algebra
$$\begin{array}{ll llll}
\mathfrak{N}^6_1 &:& e_1 e_1 = e_3 & e_1 e_2=e_4 & e_2 e_1=e_5  & e_2 e_2=e_6
\end{array}$$

\end{enumerate}

\item if    $(\mathcal{A}, \cdot)$ is a $3$-step   nilpotent algebra, then  $ 6 \leq dim(\mathcal{A}, \cdot) \leq 8$ and

\begin{enumerate}

\item if $dim(\mathcal{A}, \cdot)=6,$ then $(\mathcal{A}, \cdot)$ is isomorphic to  one of the following algebras:
\begin{longtable}{llllllll}
$\mathfrak{Z}^6_1$ &$:$& $e_1 e_2=e_3$&  $e_2e_1=e_4$& $e_2e_2=e_5$ & $e_1e_5=e_6$ &
$e_5e_1=-e_6$ \\
&& $e_2e_4=-2 e_6$ & $e_4e_2=-e_6 $& $e_2e_3=e_6$&
$e_3e_2=2e_6$\\
$\mathfrak{Z}^6_2$ &$:$& $e_1 e_2=e_3$&  $e_2e_1=e_4$& $e_2e_2=e_5+e_6$ & $e_1e_5=e_6$ &
$e_5e_1=-e_6$ \\
&& $e_2e_4=-2 e_6$ & $e_4e_2=-e_6 $& $e_2e_3=e_6$&
$e_3e_2=2e_6$
\end{longtable}

\item if $dim(\mathcal{A}, \cdot)=7,$ then $(\mathcal{A}, \cdot)$ is isomorphic to  the following algebra
\begin{longtable}{lllllllllllllll}
$\mathfrak{Z}^7_1$ &$:$& $e_1 e_1=e_7$&$e_1 e_2=e_3$&  $e_2e_1=e_4$& $e_2e_2=e_5$ & $e_1e_5=e_6$ \\
&&$e_5e_1=-e_6$ 
& $e_2e_4=-2 e_6$ & $e_4e_2=-e_6 $& $e_2e_3=e_6$&
$e_3e_2=2e_6$
\end{longtable}

\item if $dim(\mathcal{A}, \cdot)=8,$ then $(\mathcal{A}, \cdot)$ is isomorphic to  the following algebra
\begin{longtable}{lllllllllllllll}
$\mathfrak{Z}^8_1$ &$:$&
$e_1 e_1 = e_3$ & $e_1 e_2=e_4$ & $e_1 e_4 = 2 e_7$ &$e_1 e_5 = -e_7$ & $e_1 e_6 = e_8$ & $e_2 e_1=e_5$  \\
&& $e_2 e_2=e_6$ & $e_2 e_3 = -e_7$ & $e_2 e_4 = e_8$ & $e_2 e_5 =-2 e_8$ & $e_3 e_2 = e_7$ & $e_4 e_1 = e_7$ \\
&& $e_4 e_2 = 2e_8$ & $e_5 e_1 = -2e_7$ & $e_5 e_2 =- e_8$ & $e_6 e_1 = -e_8$ &

\end{longtable}

\end{enumerate}
\end{enumerate}

\end{theorem}

\dem
The first part of the theorem is following from \cite{calderon}.
The proof of the second part of the theorem will be done the usual way by calculations of symmetric Zinbiel central extensions
(it is known as the Skjelbred-Sund method 
(about the description of the method for an arbitrary variety see 
\cite{klp20} and for Zinbiel case see \cite{centr3zinb}).

Thanks to \cite{centr3zinb}, there are no non-trivial symmetric Zinbiel central extensions of algebras $\mathfrak{N}^3_i$ and $\mathfrak{N}^4_i.$
Hence, we will consider only algebras $\mathfrak{N}^5_i$ and $\mathfrak{N}^6_1.$

\begin{itemize}
    \item  The algebra $\mathfrak{N}^5_1$ has only two $1$-dimensional
    and one $2$-dimensional central extensions:

\begin{longtable}{llllllll}
$\mathfrak{Z}^6_1$ &:& $e_1 e_2=e_3$&  $e_2e_1=e_4$& $e_2e_2=e_5$ & $e_1e_5=e_6$ &
$e_5e_1=-e_6$ \\
&& $e_2e_4=-2 e_6$ & $e_4e_2=-e_6 $& $e_2e_3=e_6$&
$e_3e_2=2e_6$ \\

$\mathfrak{Z}^6_2$ &:& $e_1 e_1=e_6$& $e_1 e_2=e_3$&  $e_2e_1=e_4$& $e_2e_2=e_5$ & $e_1e_5=e_6$ \\
&&$e_5e_1=-e_6$ 
& $e_2e_4=-2 e_6$ & $e_4e_2=-e_6 $& $e_2e_3=e_6$&
$e_3e_2=2e_6$\\

$\mathfrak{Z}^7_1$ &:& $e_1 e_1=e_7$&$e_1 e_2=e_3$&  $e_2e_1=e_4$& $e_2e_2=e_5$ & $e_1e_5=e_6$ &\\
&&$e_5e_1=-e_6$ 
& $e_2e_4=-2 e_6$ & $e_4e_2=-e_6 $& $e_2e_3=e_6$&
$e_3e_2=2e_6$

\end{longtable}

\item The algebras $\mathfrak{N}^5_2,\mathfrak{N}^5_3$ and $\mathfrak{N}^5_{4,\lambda}$ have no  central extensions.

\item The algebra $\mathfrak{N}^6_1$ has only one $2$-dimensional  central extension:
\begin{longtable}{llllllll}
$\mathfrak{Z}^8_1$ &$:$&
$e_1 e_1 = e_3$ & $e_1 e_2=e_4$ & $e_1 e_4 = 2 e_7$ & $e_1 e_5 = -e_7$ & $e_1 e_6 = e_8$ & $e_2 e_1=e_5$  \\
&& $e_2 e_2=e_6$ & $e_2 e_3 = -e_7$ & $e_2 e_4 = e_8$ & $e_2 e_5 =-2 e_8$ & $e_3 e_2 = e_7$ & $e_4 e_1 = e_7$ \\
&& $e_4 e_2 = 2e_8$ & $e_5 e_1 = -2e_7$ & $e_5 e_2 =- e_8$ & $e_6 e_1 = -e_8$ &

\end{longtable} 
\hfill $\square$
\end{itemize}

\begin{proposition}\label{twood}
If $(\mathcal{A}, \cdot)$ is a symmetric Zinbiel superalgebra with two odd generators,
then $(\mathcal{A}, \cdot)$ is a $2$-step nilpotent superalgebra given in the first part of the theorem \ref{thalg2gen} with the following grading
$\mathcal{A}_{1}= \langle e_1, e_2 \rangle$ and
$\mathcal{A}_{0}= \langle e_{3}, \ldots, e_{2+i} \rangle,$ where $i=dim(\mathcal{A}_{0}).$

\end{proposition}

\dem
The description of $2$-step nilpotent symmetric Zinbiel superalgebras with two odd generators is following from the theorem  \ref{thalg2gen}.

If  $(\mathcal{A}, \cdot)$ is a $3$-step nilpotent symmetric Zinbiel superalgebra, then
the product of three elements $x,y,z$ of  $(\mathcal{A}, \cdot)$ is equal to zero if an elements from its is not a linear combination of generators.
Hence,  superalgebra $(\mathcal{A}, \cdot)$ without grading is equivalent to an algebra
with the following
identities
\[ \Omega = \{ (xy)z=x(yz-zy) \mbox{ and } x(yz)=(xy-yx)z\}.\]
By the standard Skjelbred-Sund method, it is possible to verification, that there are no non-trivial extensions of algebras from the first part of the theorem  \ref{thalg2gen} in the variety defined by identities from $\Omega.$
It is following that there are no $3$-step nilpotent symmetric Zinbiel superalgebras with two odd generators.
\hfill $\square$

\section{Mono and binary symmetric Zinbiel and Leibniz algebras} \label{sect2.5}

In this section, we will consider algebras over the complex field.

\begin{definition} \label{d31}
Let $\Omega$ be a variety of algebras defined by a family of polynomial identities,
then we say that an algebra  $\mathcal{A} \in \Omega_i$ if and only if each $i$-generated subalgebra of $\mathcal{A}$ gives an algebra from $\Omega.$ 
In particular, 
if $\mathcal{A} \in \Omega_1,$ then $\mathcal{A}$ is a mono-$\Omega$ algebra,
if $\mathcal{A} \in \Omega_2,$ then $\mathcal{A}$ is a binary-$\Omega$ algebra.

\end{definition} 

For example, let ${\rm Ass}$ be the class of associative algebras.
Then by Artin’s theorem, the class ${\rm Ass}_2$
coincides with the class of alternative algebras.
Albert’s theorem follows that the class ${\rm Ass}_1$ coincides with
the class of power-associative algebras.
It is easy to see that ${\rm  Lie}_1$ coincides with anticommutative algebras, 
i.e. they satisfy the identity $x^2=0.$
The identities of ${\rm Lie}_2$ are described by Gainov \cite{g}.
Below we will consider mono and binary symmetric Leibniz and Zinbiel algebras.

Thanks to \cite{g10,id21}, we have the description of mono and binary  left  Leibniz algebras.

\begin{theorem}
An algebra $\mathcal{A}$ is mono left Leibniz if and only if it satisfies the following identities:
\begin{center}
$x^2x = 0,$ \ $x^2x^2 = 0.$
\end{center}
An algebra $\mathcal{A}$ is binary left Leibniz if and only if it satisfies the
following identities:
\begin{center}
$x^2y=0,$ $ x(yx) = (xy)x+yx^2,$ \ $x(y(xy))=(xy)(xy) +y(x(xy)).$
\end{center}

\end{theorem}

Thanks to \cite{ims21}, we have the description of mono and binary  left  Zinbiel algebras.

\begin{theorem}
An algebra $\mathcal{A}$ is mono left Zinbiel if and only if it satisfies the following identities:
\begin{center}
$xx^2 = 2x^2x,$ \ $x^2x^2 = 3(x^2x)x.$
\end{center}
An algebra $\mathcal{A}$ is binary left Zinbiel if and only if it satisfies the
following identities:
\begin{center}
$x(yx) = (xy+yx)x,$ \ $x(xy) = 2x^2y.$
\end{center}

\end{theorem}

As some trivial corollaries from the previous theorems, we have the following propositions.

\begin{proposition} 
The variety of mono  symmetric Zinbiel algebras coincides with the variety of mono symmetric Leibniz algebras and it is defined by the following identities 
\begin{center}
    $x^2x=0,$ \ $xx^2=0, $ \ $x^2x^2=0.$    
\end{center}
\end{proposition}

\begin{proposition}\label{21}
The variety of binary  symmetric Zinbiel algebras is defined by the following identities
\begin{center}
$x(yx) = (xy+yx)x,$ \ $x(xy) = 2x^2y,$\\
$(xy)x = x(yx+xy),$ \ $(yx)x = 2yx^2.$

\end{center}\end{proposition} 

\begin{proposition}\label{22}
The variety of binary  symmetric Zinbiel algebras is defined by the following identities
\begin{center}
$x^2y=0,$ \ $x(yx) = (xy)x,$ \ $x(y(xy))=(xy)(xy)+y(x(xy))$\\
$yx^2 = 0,$ \ $(xy)(xy) = ((xy)x)y+x((xy)y).$

\end{center}\end{proposition}

\begin{proposition} 
The intersection of varieties  of binary   symmetric Zinbiel and Leibniz  algebras  
it is defined by the following identities 
\begin{center}
    $(x_1x_2)x_3=(-1)^{\sigma}(x_{\sigma(1)}x_{\sigma(2)})x_{\sigma(3)}$ and 
    $x_1(x_2x_3)=(-1)^{\sigma}x_{\sigma(1)}(x_{\sigma(2)}x_{\sigma(3)}),$ \ $\sigma \in \mathbb{S}_3.$
\end{center} 
\end{proposition} 

\dem
 Thanks to two previous  propositions, we have that each algebra from the intersection of
symmetric Zinbiel and Leibniz  algebras  satisfies the following 
identities 
\begin{center}
$x^2y=0,$ $x(xy)=0,$  \ $yx^2=0,$ $(yx)x=0.$
\end{center}
Hence, by linearization, it satisfies 
\begin{center}
$(x_1x_2)x_3=-(x_2x_1)x_3,$ $x_1(x_2x_3)=-x_2(x_1x_3),$  \ $x_1(x_2x_3)=-x_1(x_3x_2),$ \  $(x_1x_2)x_3=-(x_1x_3)x_2.$
\end{center}
It is easy to see that each algebra that satisfies the previous identities also satisfies all identities from two previous propositions.
Hence, we complete the proof of the statement. \hfill $\square$

\begin{theorem}
Let ${\rm SZ}$ and ${\rm SL}$ be the varieties of symmetric Zinbiel and symmetric Leibniz algebras.
Then, the following inclusions are strict.
$$
\begin{array}{ccccccc} 

            & & {\rm SZ}  & \subset &   {\rm SZ}_2       \\
           & \rotatebox{45}{$\subset$}  &  &      & \rotatebox{90}{$\subset$}       & \rotatebox{-45}{$\subset$}                          \\
{\rm SL} \cap {\rm SZ}={\rm SL} \cap {\rm SZ}_2 = {\rm SZ} \cap {\rm SL}_2 & \subset         & {\rm Ass } \cap {\rm  Lie}_1&     \subset      & {\rm SL}_2 \cap {\rm SZ}_2       & \subset  & {\rm SZ}_1={\rm SL}_1                                                             \\
           & \rotatebox{-45}{$\subset$} &     &    & \rotatebox{-90}{$\subset$}         & \rotatebox{45}{$\subset$}                      \\
           &                            & {\rm SL} & \subset & {\rm SL}_2 
\end{array}
$$
\end{theorem}

\dem
Let us prove that ${\rm SL} \cap {\rm SZ}_2$ is defined by $(xy)z=0$ and $x(yz)=0.$
Obviously, that each algebra from ${\rm SL} \cap {\rm SZ}_2$ satisfies identities from the previous proposition and two Leibniz identities. Hence,
\begin{center}
    $2(xy)z=(xy)z+(xz)y+x(yz)=x(yz)$ and 
    $2x(yz)=x(yz)+(xy)z+y(xz)=(xy)z,$
\end{center}
which gives $(xy)z=0$ and $x(yz)=0.$

${\rm SZ} \cap {\rm SL}_2$ is also defined by $(xy)z=0$ and $x(yz)=0,$ because 
each algebra from ${\rm SZ} \cap {\rm SL}_2$ satisfies identities from the previous proposition.
Hence, it satisfies
\begin{center}
    $(xy)z=x(yz+zy)=0$ and 
    $x(yz)=(xy+yx)z=0.$
\end{center}

Obviously, we have the equalities ${\rm SL} \cap {\rm SZ}={\rm SL} \cap {\rm SZ}_2 = {\rm SZ} \cap {\rm SL}_2.$

All indicated inclusions are following from previous propositions and other known results.
To prove the strict inclusion $\omega \subset \Omega$ of varieties of algebras, 
we will indicate an algebra $\mathcal{A} \in \Omega \setminus \omega.$  

\begin{longtable}{|rcl| llll|}
\hline
${\rm SL} \cap {\rm SZ}$&$ \subset $&${\rm SZ}$  &  
$e_1 e_2=e_3$&  $e_1e_5=e_6$ & $e_2e_1=e_4$&\\
${\rm SL}_2 \cap {\rm SZ}_2$&$ \subset$&$  {\rm SZ}_2$&   $e_2e_2=e_5$ &    $e_2e_3=e_6$&  $e_2e_4=-2 e_6$ &\\ 
${\rm SL}_2$&$  \subset$&$ {\rm SL}_1$ &      $e_3e_2=2e_6$& $e_4e_2=-e_6 $& $e_5e_1=-e_6$ &

\\

\hline
${\rm SL} \cap {\rm SZ}$&$ \subset$&$ {\rm SL}$   &   
$e_1 e_2=e_3$&  $e_2e_1=-e_3$& $e_2e_3=e_4$ & $e_3e_2=-e_4$
\\

\hline

${\rm SL} \cap {\rm SZ}$&$ \subset$&${\rm Ass } \cap {\rm  Lie}_1$  &  
$e_1 e_2=e_4$& $e_1e_3=e_5$ & $e_1e_6=e_7$&  $e_2e_1=-e_4$ \\
&&& $e_2 e_3=e_6$   & $e_2 e_5=-e_7$  & $e_3e_1=-e_5$ &  $e_3e_2=-e_6$   \\
&&& 
  $e_3e_4=e_7$ & $e_4e_3=e_7$ &  $e_5e_2=-e_7$ & $e_6e_1=e_7$\\

\hline

${\rm Ass } \cap {\rm  Lie}_1$&$ \subset$&$ {\rm SL}_2 \cap {\rm SZ}_2$  &   
$e_1 e_2=e_4$ & $e_1e_3=e_5$ & $e_1e_6=e_7$&  \\
${\rm SL}$&$  \subset$&$ {\rm SL}_2$ &    $e_2e_1=-e_4$ &
  $e_2 e_3=e_6$&  $e_2 e_5=-e_7$ & \\
${\rm SZ}$&$  \subset$&$ {\rm SZ}_2$ &  $e_3e_1=-e_5$ &  $e_3e_2=-e_6$  &   $e_3e_4=e_7$  & \\

\hline

${\rm SZ}_2$&$  \subset$&$ {\rm SZ}_1$  &  
$e_1e_2=e_2$ &  &&\\
${\rm SL}_2 \cap {\rm SZ}_2$&$ \subset$&$ {\rm SL}_2$ & $e_2e_1=-e_2$  &&&   \\
\hline
\end{longtable}

\hfill $\square$

From the proof of the present theorem
it is easy to see that the variety ${\rm SL}_2 \cap {\rm SZ}_2$ does not contain ${\rm SL}$ and ${\rm SZ}$ as subvarieties.
It gives the following interesting open question.

\begin{oquestion}
Is  ${\rm SL}_1$ the minimal variety that contains the varieties 
${\rm SL}$ and ${\rm SZ}$ (resp. ${\rm SL}_2$ and ${\rm SZ}_2$) as subvarieties?

\end{oquestion}

\section{Quadratic Zinbiel superalgebras} \label{sect3}
\begin{definition} \label{d31}
A triple $(\mathcal{A}, \cdot , \mathfrak{B})$ is said to be a quadratic Zinbiel superalgebra if $(\mathcal{A}, \cdot )$ is a left (or a right) Zinbiel superalgebra endowed with a non-degenerate bilinear form $\mathfrak{B}$ satisfying the following properties:
\begin{enumerate}
\item[\bf (i)]
$\mathfrak{B}(x,y)=(-1)^{\vert x\vert \vert y\vert} \mathfrak{B}(y,x),~\forall x\in \mathcal{A}_{\vert x\vert},~y\in \mathcal{A}_{\vert y\vert}$, i.e. $\mathfrak{B}$ is supersymmetric;
\item[\bf (ii)] $\mathfrak{B}(x  y,z)=\mathfrak{B}(x,y  z),~\forall x,y,z\in \mathcal{A}$, i.e. $\mathfrak{B}$ is invariant;
\item[\bf (iii)] $\mathfrak{B}(\mathcal{A}_{\bar{0}},\mathcal{A}_{\bar{1}})=\mathfrak{B}(\mathcal{A}_{\bar{1}},\mathcal{A}_{\bar{0}})=\{0\}$, i.e. $\mathfrak{B}$ is even.
\end{enumerate}
In this case, $\mathfrak{B}$ is called an invariant scalar product on $\mathcal{A}$.
\end{definition}
\begin{proposition} \label{QUSY}
Let $(\mathcal{A},\cdot )$ be a left (resp. right) Zinbiel superalgebra. If $\mathcal{A}$ is quadratic, then $\mathcal{A}$ is symmetric.
\end{proposition}
\dem
Suppose that $(\mathcal{A},\cdot)$ is a left Zinbiel superalgebra. Then, for all
$x,  y,  z, t\in \mathcal{A}_{\bar{0}} \cup\mathcal{A}_{\bar{1}},$
we have
$$0=\mathfrak{B}((x y) z-x  (y z)-(-1)^{\vert y\vert \vert z\vert} x  (z  y),t)=
\mathfrak{B}(x,y  (z t)-(y  z)  t-(-1)^{\vert y\vert \vert z\vert}(z  y)  t).$$
The fact that $\mathfrak{B}$ is non-degenerate implies that
$y (z  t)-(y  z)  t-(-1)^{\vert y\vert \vert z\vert}(z  y)  t=0$. 
Consequently, $\mathcal{A}$ is a right Zinbiel superalgebra. By the same way, we prove the result for right Zinbiel superalgebra. \hfill $\square$
\begin{proposition} \label{p38}
Let $\mathcal{A}=\mathcal{A}_{\bar{0}}\oplus \mathcal{A}_{\bar{1}}$ be a symmetric Zinbiel superalgebra. Then, $\mathcal{A}$ is quadratic if and only if the adjoint and the co-adjoint representations of $\mathcal{A}$ are equivalent and $dim(\mathcal{A}_{\bar{1}})$ is even.
\end{proposition}
\dem
The proof is similar to that in \cite{HS04}.~~\hfill $\square$\\

The following result comes from Proposition \ref{corep}.
\begin{proposition}\label{quadra}
Let $(\mathcal{A}, \cdot )$ be a left (resp. right) Zinbiel superalgebra. If $\mathcal{A}$ is quadratic, then $\mathcal{A}$ is $2$-step nilpotent.
\end{proposition}

\section{Extensions of quadratic Zinbiel superalgebras} \label{sect3}
In this section, we introduce the notion of a double extension of quadratic Zinbiel superalgebras by the one-dimensional $\mathbb{Z}_{2}$-graded vector space, which is performed in two distinct steps, a central extension followed by a semi-direct product. In the first step, we are going to introduce the notion of central extension by the one-dimensional $\mathbb{Z}_{2}$-graded vector space.

\begin{proposition}\label{oprt}
Let $(\mathcal{A}, \cdot , \mathfrak{B})$ be a quadratic Zinbiel superalgebra and $\Omega:\mathcal{A}\times \mathcal{A}\rightarrow \mathbb{K}$ be a bilinear map of degree $\alpha$. Then, there exists a homogeneous endomorphism $\delta$ of $\mathcal{A}$ of degree $\alpha$ such that $\Omega(x,y)=\mathfrak{B}(\delta(x),y),~\forall x,y\in \mathcal{A}$. The map $\Omega$ is a homogeneous scalar Zinbiel 2-cocycle of $\mathcal{A}$ of degree $\alpha$ if and only if, for all $x,y\in \mathcal{A}_{\vert x\vert}\times \mathcal{A}_{\vert y\vert}$, the following assertions hold:
$${\bf (i)}~\delta(x  y)=\delta(x)  y+(-1)^{\vert y\vert (\vert x\vert+\alpha)}y \delta(x)\quad 
and \quad
{\bf (ii)}~\delta(x  y)=-(-1)^{\alpha\vert x\vert}x  \delta(y).$$
\end{proposition}
\dem ${\bf (i)}$ It is easy to see the following
$$\Omega(x,y  z)+(-1)^{\vert y\vert \vert z\vert}\Omega(x,z  y)-\Omega(x  y,z)=\mathfrak{B}(\delta(x)  y+(-1)^{\vert y\vert (\vert x\vert+\alpha)}y \delta(x)-\delta(x y),z).$$ 
\\${\bf (ii)}$ Note that $~\forall x,y,z\in \mathcal{A}_{\vert x\vert}\times \mathcal{A}_{\vert y\vert}\times \mathcal{A}_{\vert z\vert}$ we have
$$\Omega(x,y z)+(-1)^{\vert z\vert (\vert x\vert + \vert y\vert)} \Omega(z x,y)=
(-1)^{\vert x\vert (\vert y\vert + \vert z\vert)+\alpha \vert y\vert}\mathfrak{B}(y,(-1)^{\alpha \vert z\vert}z \delta(x)+\delta(z  x)).$$ \hfill $\square$

\begin{proposition}\label{extcent}
Let $(\mathcal{A}, \cdot , \mathfrak{B})$ be a quadratic Zinbiel superalgebra, $\mathcal{V}=\mathbb{K}d$ be the one-dimensional $\mathbb{Z}_{2}$-graded vector space and $\delta$ be a homogeneous endomorphism of $\mathcal{A}$ of degree $\alpha$. Define the  bilinear map $\Omega:\mathcal{A}\times \mathcal{A}\rightarrow \mathcal{V^{*}}$ by $\Omega(x,y)=\mathfrak{B}(\delta(x),y)d^{*}$. Then, the $\mathbb{Z}_{2}$-graded vector space $\widetilde{\mathcal{A}}=\mathcal{A}\oplus \mathcal{V^{*}}$ endowed with the product defined by
$$(x+\lambda d^{*})\star_{\Omega} (y+\lambda' d^{*})=x  y+\Omega(x,y),~~\mbox{for~all}~(x+\lambda d^{*}),(y+\lambda' d^{*})\in \widetilde{\mathcal{A}},$$ is a symmetric Zinbiel superalgebra if and only if $\Omega$ is a homogeneous Zinbiel 2-cocycle of $\mathcal{A}$ on the trivial $\mathcal{A}-$module $\mathcal{V^{*}}$ of degree $\alpha$. It is termed the central extension of $\mathcal{A}$ by $\mathcal{V^{*}}$ by means of $\Omega$.
\end{proposition}
\dem We only need calculation, where we use Proposition \ref{oprt}. \hfill $\square$\\

We will define, in the following, another type of extension of symmetric Zinbiel superalgebras.
Let us consider $(\mathcal{A},\cdot)$ a symmetric Zinbiel superalgebra and $\mathcal{V}=\mathbb{K}d$ a one-dimensional $\mathbb{Z}_{2}$-graded vector space. Let $\delta,D$ be two homogeneous endomorphisms of $\mathcal{A}$ of degree $\vert d\vert$ and $\mathfrak{a}_{0}\in {\rm Ann}(\mathcal{A})\cap \mathcal{A}_{\overline{0}}$. Consider the $\mathbb{Z}_{2}$-graded vector space $\overline{\mathcal{A}}=\mathcal{A}\oplus \mathcal{V}$ on which we define the following product:
\begin{equation}\label{equasedipr}
d\mydot d=\mathfrak{a}_{0};\quad x\mydot y=x y;\quad
d\mydot x=\delta(x);\quad x\mydot d=D(x),~~\forall x,y\in \mathcal{A}.
\end{equation}
\begin{proposition}
The $\mathbb{Z}_{2}$-graded vector space $\overline{\mathcal{A}}$ endowed with the product (\ref{equasedipr}) is a symmetric Zinbiel superalgebra if and only if, for all $x,y\in \mathcal{A}_{\vert x\vert}\times \mathcal{A}_{\vert y\vert}$, the following conditions hold:
\begin{eqnarray*}
&~&\delta(x y)=(-1)^{\vert d\vert \vert x\vert}x  \delta(y)=\delta(x) y+(-1)^{\vert d\vert \vert x\vert}D(x)  y=-(-1)^{\vert y\vert (\vert x\vert+\vert d\vert)}D(y) x;\\
&~&\delta(x) y=(-1)^{\vert x\vert \vert y\vert}\delta(y) x=-(-1)^{\vert d\vert (\vert x\vert+\vert y\vert)}x D(y);\\
&~&D(x y)=(-1)^{\vert d\vert \vert y\vert}D(x) y;\\
&~&\delta^{2}(x)=(-1)^{\vert d\vert}\delta^{2}(x)=-D^{2}(x)=(1+(-1)^{\vert d\vert}) \mathfrak{a}_{0} x;\\
&~&\delta\circ D(x)=-(-1)^{\vert d\vert}D\circ \delta(x)=(-1)^{\vert d\vert \vert x\vert}x \mathfrak{a}_{0}=
-(-1)^{\vert d\vert (\vert x\vert+\vert d\vert)}\mathfrak{a}_{0} x;~~\delta(\mathfrak{a}_{0})=D(\mathfrak{a}_{0})=0.
\end{eqnarray*}
In this case, $(\delta,D,\mathfrak{a}_{0})$ is called "an admissible triple" and the symmetric Zinbiel superalgebra $\overline{\mathcal{A}}$ is termed the semi-direct product of $\mathcal{A}$ by $\mathcal{V}$ by means of $(\delta,D,\mathfrak{a}_{0})$.
\end{proposition}
\dem The proof is a straightforward computation. \hfill $\square$\\

Now, we are in a position to introduce the double extensions of quadratic Zinbiel superalgebras. We start by presenting
the concept of even double extension of these superalgebras by the one-dimensional
Lie algebra.
\begin{theorem} \label{EvenDoubExtZinb}
Let $(\mathcal{A},\cdot , \mathfrak{B})$ be a quadratic Zinbiel superalgebra, $\mathcal{V}=\mathbb{K}d$ be the one-dimensional vector space and $(\delta,\delta^{*},\mathfrak{a}_{0})$ be an admissible triple of $\mathcal{A}$ such that
$\delta(\mathcal{A})\subseteq {\rm Ann}(\mathcal{A}),~\delta(\mathcal{A}^{2})=\{0\},~
\delta^{2}=\delta\circ \delta^{*}=\delta^{*}\circ \delta=\{0\},~\mathfrak{a}_{0}\in {\rm Ann}(\mathcal{A})\cap \mathcal{A}_{\bar{0}}$ and $\mathfrak{B}(\mathfrak{a}_{0},\mathfrak{a}_{0})=0$, where $\delta^{*}$ is the adjoint of $\delta$ with respect to $\mathfrak{B}$. Then, the $\mathbb{Z}_{2}$-graded vector space $\overline{\mathcal{A}}=\mathcal{V}^{*}\oplus \mathcal{A}\oplus \mathcal{V}$ endowed with the following product:
\begin{longtable}{lll}
$d\mydot d=\mathfrak{a}_{0}+\alpha d^{*};$ & $d\mydot x=\delta(x)+\mathfrak{B}(x,\mathfrak{a}_{0})d^{*};$ \\
$x\mydot y=xy+\mathfrak{B}(\delta(x),y)d^{*};$  & $x\mydot d=\delta^{*}(x)+\mathfrak{B}(x,\mathfrak{a}_{0})d^{*},$ & $\forall x,y\in \mathcal{A},~\alpha\in \mathbb{K},$
\end{longtable}
is a symmetric Zinbiel superalgebra. Moreover, the even bilinear form
$\overline{\mathfrak{B}}:\overline{\mathcal{A}}\times \overline{\mathcal{A}}\rightarrow \mathbb{K}$, defined by $\overline{\mathfrak{B}}|_{\mathcal{A}\times\mathcal{A}}=\mathfrak{B};~~\overline{\mathfrak{B}}(d,d^{*})=
\overline{\mathfrak{B}}(d^{*},d)=1,$
is an associative scalar product on $\overline{\mathcal{A}}$.

The Zinbiel superalgebra $(\overline{\mathcal{A}},\overline{\mathfrak{B}})$ is called the
even double extension of $\mathcal{A}$ by $\mathcal{V}$ by means of the admissible
triple $(\delta,\delta^{*},\mathfrak{a}_{0})$.
\end{theorem}
\dem Since $\delta(\mathcal{A}^{2})=\{0\}$ and $\delta(\mathcal{A})\subseteq {\rm Ann}(\mathcal{A})$, then Proposition \ref{oprt} implies that the even map $\Omega:\mathcal{A}\times \mathcal{A}\rightarrow \mathcal{V}^{*},~(x,y)\longmapsto \Omega(x,y)=\mathfrak{B}(\delta(x),y)d^{*}$, is an even Zinbiel 2-cocycle of $\mathcal{A}$ on the trivial $\mathcal{A}-$module $\mathcal{V}^{*}$. So, $\widetilde{\mathcal{A}}=\mathcal{A}\oplus \mathcal{V}^{*}$ endowed with the product defined by:$$(x+f)\star_{\Omega}(y+g)=x y+\mathfrak{B}(\delta(x),y)d^{*},~\forall x,y\in \mathcal{A},~f,g\in \mathcal{V}^{*},$$ is a symmetric Zinbiel superalgebra, central extension of $\mathcal{A}$ by means of $\Omega$. Now, we define two even endomorphisms $\delta_{1}$ and $\delta_{2}$ of $\widetilde{\mathcal{A}}$, respectively, by 
\begin{center}$\delta_{1}(x+\lambda d^{*})=\delta(x)+\mathfrak{B}(x,\mathfrak{a}_{0})d^{*}$ and $\delta_{2}(x+\lambda d^{*})=\delta^{*}(x)+\mathfrak{B}(x,\mathfrak{a}_{0})d^{*},~\forall x\in \mathcal{A},~\lambda \in \mathbb{K}.$\end{center}
Let us consider the element $\mathfrak{a}_{1}=\mathfrak{a}_{0}+\alpha d^{*}\in \widetilde{\mathcal{A}}$, where $\alpha$ is a fixed scalar in $\mathbb{K}$. So, $\mathfrak{a}_{1}\in {\rm Ann}(\widetilde{\mathcal{A}})\cap \widetilde{\mathcal{A}}_{\bar{0}}$.
Since 
\begin{center}$\delta(\mathcal{A})\subseteq {\rm Ann}(\mathcal{A}),~\delta(\mathcal{A}^{2})=\{0\},~
\delta^{2}=\delta\circ \delta^{*}=\delta^{*}\circ \delta=\{0\},~\mathfrak{a}_{0}\in {\rm Ann}(\mathcal{A})\cap \mathcal{A}_{\bar{0}},~\mathfrak{B}(\mathfrak{a}_{0},\mathfrak{a}_{0})=0$,\end{center} then a simple computation proves that the triple $(\delta_{1},\delta_{2},\mathfrak{a}_{1})$ is an admissible triple of $\widetilde{\mathcal{A}}$. Consequently, we can consider the symmetric Zinbiel superalgebra $\overline{\mathcal{A}}=\widetilde{\mathcal{A}}\oplus \mathcal{V}=\mathcal{V}^{*}\oplus \mathcal{A}\oplus \mathcal{V}$, semi-direct product of $\widetilde{\mathcal{A}}$ by $\mathcal{V}$. In addition, it is easy to show that the bilinear form $\overline{\mathfrak{B}}:\overline{\mathcal{A}}\times \overline{\mathcal{A}}\rightarrow \mathbb{K}$, defined by $\overline{\mathfrak{B}}|_{\mathcal{A}\times\mathcal{A}}=\mathfrak{B};~\overline{\mathfrak{B}}(d,d^{*})=1$, is an associative scalar product on $\overline{\mathcal{A}}$. \hfill $\square$\\


Our next goal is to establish the converse of Theorem \ref{EvenDoubExtZinb}.
\begin{theorem} \label{EvenIndDescrZinb}
Let $(\mathcal{A},\mydot,\mathfrak{B})$ be a quadratic Zinbiel superalgebra of dimension $(n,m)$ such that $n\geq 2$. If ${\rm Ann}(\mathcal{A})\cap \mathcal{A}_{\bar{0}}\neq\{0\}$,
then $\mathcal{A}$ is isomorphic to an even double extension of a quadratic Zinbiel superalgebra $(\mathcal{H},\cdot,\mathfrak{B}_{\mathcal{H}})$ of dimension $(n-2,m)$ by the one-dimensional vector space.
\end{theorem}
\dem Suppose that ${\rm Ann}(\mathcal{A})\cap \mathcal{A}_{\bar{0}}\neq\{0\}$, then there exists a non-zero element
$e$ of ${\rm Ann}(\mathcal{A})\cap \mathcal{A}_{\bar{0}}$. Denote by
$\mathcal{J}=\mathbb{K}e$ and by $\mathcal{J}^{\bot}$ the orthogonal of
$\mathcal{J}$ with respect to $\mathfrak{B}$. Then, $\mathcal{J}\subseteq \mathcal{J}^{\bot}$. Since
$\mathfrak{B}$ is even and non-degenerate, then there exists $d\in \mathcal{A}_{\bar{0}}\backslash \{0\}$ such that $\mathfrak{B}(e,d)=\mathfrak{B}(d,e)=1$. Let $\mathcal{V}=\mathbb{K}d$ and $\mathcal{H}=(\mathcal{J}\oplus \mathcal{V})^{\bot}$, then
$\mathcal{A}=\mathcal{J}\oplus \mathcal{H}\oplus \mathcal{V}$ and $\mathcal{J}^{\bot}=\mathcal{J}\oplus \mathcal{H}$
is an ideal of $\mathcal{A}$. Let
$x,y\in \mathcal{H}_{\vert x\vert}\times\mathcal{H}_{\vert y\vert}$,
then $x\mydot y=\Phi(x,y)e+x y$, where
$\Phi:\mathcal{H}\times \mathcal{H}\rightarrow \mathbb{K}$ is an even bilinear
form and $``\cdot":\mathcal{H}\times \mathcal{H}\rightarrow \mathcal{H}$ is an
even bilinear map. It is easy to show that $(\mathcal{H},\cdot)$ is a
symmetric Zinbiel superalgebra and that the even bilinear form
$\mathfrak{B}_{\mathcal{H}}=\mathfrak{B}|_{\mathcal{H}\times \mathcal{H}}$ is an associative scalar product
on $(\mathcal{H},\cdot)$. Therefore,
the even bilinear form $\Phi$ is an element of $(Z^2_{\rm SZinb}(\mathcal{H},\mathbb{K}))_{\bar{0}}$. The fact that $\mathcal{J}^{\bot}=\mathcal{J}\oplus \mathcal{H}$ is an ideal of $\mathcal{A}$ implies that
\begin{eqnarray*}
&~& d\mydot d=\alpha e+\mathfrak{a}_{0}+\lambda d,~\hbox{where}~\alpha,\lambda \in \mathbb{K},~
\mathfrak{a}_{0}\in \mathcal{H}_{\bar{0}},\\
&~& x\mydot d=\Psi(x)e+D(x),~\hbox{where}~D\in (End(\mathcal{H}))_{\bar{0}},~
\Psi\in (\mathcal{H}^{*})_{\bar{0}},\\
&~& d\mydot x=\Theta(x)e+\delta(x),~\hbox{where}~\delta\in (End(\mathcal{H}))_{\bar{0}},~
\Theta\in (\mathcal{H}^{*})_{\bar{0}}.
\end{eqnarray*}
Let $x,y\in \mathcal{H}_{\vert x\vert}\times\mathcal{H}_{\vert y\vert}$. Since $\mathfrak{B}$ is supersymmetric and associative, then $\mathfrak{B}(x\mydot y,d)=\mathfrak{B}(x,y\mydot d)=\mathfrak{B}(d\mydot x,y)$. It follows that $\Phi(x,y)=\mathfrak{B}(x,D(y))=\mathfrak{B}(\delta(x),y)$. So, $D=\delta^{*}$. Moreover, $\Psi(x)=\mathfrak{B}(x\mydot d,d)=\mathfrak{B}(d,d\mydot x)=\Theta(x)=\mathfrak{B}(x,\mathfrak{a}_{0})$ and $\lambda=\mathfrak{B}(d\mydot d,e)=\mathfrak{B}(d,d\mydot e)=0$. Then, $d\mydot d=\alpha e+\mathfrak{a}_{0}$. Since $\mathcal{A}$ is a left Zinbiel superalgebra satisfying:
$$\mathscr{X}\mydot(\mathscr{Y}\mydot \mathscr{Z})+(-1)^{\vert z\vert (\vert x\vert + \vert y\vert)}(\mathscr{Z}\mydot \mathscr{X})\mydot \mathscr{Y}=0,~~
\mbox{for~all}~\mathscr{X}\in \mathcal{A}_{{\vert x\vert}},
\mathscr{Y}\in \mathcal{A}_{{\vert y\vert}},\mathscr{Z} \in \mathcal{A}_{{\vert z\vert}},$$
then a simple computation proves that $(\delta,D,\mathfrak{a}_{0})$ is an admissible
triple of $\mathcal{H}$, 
\begin{center}$\delta(\mathcal{H})\subseteq {\rm Ann}(\mathcal{H}),~\delta(\mathcal{H}^{2})=\{0\},~
\delta^{2}=\delta\circ D=D\circ \delta=\{0\},~\mathfrak{a}_{0}\in {\rm Ann}(\mathcal{H})\cap \mathcal{H}_{\bar{0}}$ and $\mathfrak{B}(\mathfrak{a}_{0},\mathfrak{a}_{0})=0$.
\end{center}Therefore, we can consider the quadratic Zinbiel superalgebra $\overline{\mathcal{A}}=\mathcal{V}^{*}\oplus \mathcal{H}\oplus \mathcal{V}$, even double extension of $(\mathcal{H},\cdot,\mathfrak{B}_{\mathcal{H}})$ by $\mathcal{V}$ by means of $(\delta,D,\mathfrak{a}_{0})$.

It is clear that the map 
\begin{center}$\Gamma:\mathcal{A}\rightarrow \overline{\mathcal{A}},~\lambda e+x+\lambda'd\longmapsto \lambda d^{*}+x+\lambda'd$,
\end{center} is an isomorphism of Zinbiel superalgebras. Then, $\mathcal{A}$ is isomorphic to the even double extension of $\mathcal{H}$ by $\mathcal{V}$. \hfill $\square$\\

In the second part of this section, we will present the notion of an odd double
extension of quadratic Zinbiel superalgebras by the one-dimensional $\mathbb{Z}_{2}$-graded vector space with even part zero.
\begin{theorem} \label{OddDoubExtZinb}
Let $(\mathcal{A},\cdot,\mathfrak{B})$ be a quadratic Zinbiel superalgebra, $\mathcal{N}=\mathcal{N}_{\bar{1}}=\mathbb{K}d$ be
the one-dimensional $\mathbb{Z}_{2}$-graded vector space with even part zero and $(\delta,D,\mathfrak{a}_{0})$ be an admissible triple of $\mathcal{A}$ such that \begin{center}$\delta(\mathcal{A})\subseteq {\rm Ann}(\mathcal{A}),~\delta(\mathcal{A}^{2})=\{0\},~
\delta\circ D=D\circ \delta=\{0\},~\mathfrak{a}_{0}\in {\rm Ann}(\mathcal{A})\cap \mathcal{A}_{\bar{0}}$ and $\mathfrak{B}(\mathfrak{a}_{0},\mathfrak{a}_{0})=0$, 
\end{center} where $D\in (End(\mathcal{A}))_{\bar{1}}$ verifying $\mathfrak{B}(\delta(x),y)=(-1)^{\vert x\vert + \vert y\vert} \mathfrak{B}(x,D(y))$
for all $x,y\in \mathcal{A}_{\vert x\vert}\times \mathcal{A}_{\vert y\vert}$.
Then, the $\mathbb{Z}_{2}$-graded vector space $\overline{\mathcal{A}}=\mathcal{N}^{*}\oplus \mathcal{A}\oplus \mathcal{N}$ endowed with the following product:
\begin{longtable}{lll}
$d\mydot d=\mathfrak{a}_{0};$ & $d\mydot x=\delta(x)-\mathfrak{B}(x,\mathfrak{a}_{0})d^{*};$  \\
 $x\mydot y=xy-\mathfrak{B}(\delta(x),y)d^{*};$& $x\mydot d=D(x)+\mathfrak{B}(x,\mathfrak{a}_{0})d^{*},$ & $\forall x,y\in \mathcal{A},$
\end{longtable}
is a symmetric Zinbiel superalgebra. Moreover, the even bilinear form
$\overline{\mathfrak{B}}:\overline{\mathcal{A}}\times \overline{\mathcal{A}}\rightarrow \mathbb{K}$, defined by $\overline{\mathfrak{B}}|_{\mathcal{A}\times\mathcal{A}}=\mathfrak{B};~~\overline{\mathfrak{B}}(d^{*},d)=
-\overline{\mathfrak{B}}(d,d^{*})=1$,
is an associative scalar product on $\overline{\mathcal{A}}$.

The Zinbiel superalgebra $(\overline{\mathcal{A}},\overline{\mathfrak{B}})$ is termed the
odd double extension of $\mathcal{A}$ by $\mathcal{N}$ by means of the admissible
triple $(\delta,D,\mathfrak{a}_{0})$.
\end{theorem}
\dem The proof is similar to that of Theorem \ref{EvenDoubExtZinb}. ~~\hfill $\square$\\

The following theorem is the converse of Theorem \ref{OddDoubExtZinb}.
\begin{theorem} \label{EvenIndDescrZinb}
Let $(\mathcal{A},\mydot,\mathfrak{B})$ be a quadratic Zinbiel superalgebra of dimension $(n,m)$ such that $m\geq 2$. If ${\rm Ann}(\mathcal{A})\cap \mathcal{A}_{\bar{1}}\neq\{0\}$,
then $\mathcal{A}$ is isomorphic to an odd double extension of a quadratic Zinbiel superalgebra $(\mathcal{H},\cdot,\mathfrak{B}_{\mathcal{H}})$ of dimension $(n,m-2)$ by the one-dimensional Lie
superalgebra with even part zero.
\end{theorem}
\dem Suppose that ${\rm Ann}(\mathcal{A})\cap \mathcal{A}_{\bar{1}}\neq\{0\}$, then there exists a non-zero element
$e$ of ${\rm Ann}(\mathcal{A})\cap \mathcal{A}_{\bar{1}}$. Denote by
$\mathcal{J}=\mathbb{K}e$ and by $\mathcal{J}^{\bot}$ the orthogonal of
$\mathcal{J}$ with respect to $\mathfrak{B}$. Then, $\mathcal{J}\subseteq \mathcal{J}^{\bot}$. Since
$\mathfrak{B}$ is even and non-degenerate, then there exists $d\in \mathcal{A}_{\bar{1}}\backslash \{0\}$ such that $\mathfrak{B}(e,d)=-\mathfrak{B}(d,e)=1$. Let $\mathcal{N}=\mathbb{K}d$ and $\mathcal{H}=(\mathcal{J}\oplus \mathcal{N})^{\bot}$, then
$\mathcal{A}=\mathcal{J}\oplus \mathcal{H}\oplus \mathcal{N}$ and $\mathcal{J}^{\bot}=\mathcal{J}\oplus \mathcal{H}$ is an ideal of $\mathcal{A}$. Let
$x,y\in \mathcal{H}_{\vert x\vert}\times\mathcal{H}_{\vert y\vert}$,
then $x\mydot y=\Phi(x,y)e+x y$, where $\Phi:\mathcal{H}\times \mathcal{H}\rightarrow \mathbb{K}$ is an odd bilinear
form and $``\cdot":\mathcal{H}\times \mathcal{H}\rightarrow \mathcal{H}$ is an
even bilinear map. It is easy to show that $(\mathcal{H},\cdot)$ is a
symmetric Zinbiel superalgebra and that the even bilinear form
$\mathfrak{B}_{\mathcal{H}}=\mathfrak{B}|_{\mathcal{H}\times \mathcal{H}}$ is an associative scalar product
on $(\mathcal{H},\cdot)$. Therefore,
the odd bilinear form $\Phi$ is an element of $(Z^2_{\rm SZinb}(\mathcal{H},\mathbb{K}))_{\bar{1}}$. Since $\mathcal{J}^{\bot}=\mathcal{J}\oplus \mathcal{H}$ is an ideal of $\mathcal{A}$, then
\begin{eqnarray*}
&~& d\mydot d=\mathfrak{a}_{0},~\hbox{where}~\mathfrak{a}_{0}\in \mathcal{H}_{\bar{0}},\\
&~& x\mydot d=\Psi(x)e+D(x),~\hbox{where}~D\in (End(\mathcal{H}))_{\bar{1}},~
\Psi\in (\mathcal{H}^{*})_{\bar{0}},\\
&~& d\mydot x=\Theta(x)e+\delta(x),~\hbox{where}~\delta\in (End(\mathcal{H}))_{\bar{1}},~
\Theta\in (\mathcal{H}^{*})_{\bar{0}}.
\end{eqnarray*}
Let $x,y\in \mathcal{H}_{\vert x\vert}\times\mathcal{H}_{\vert y\vert}$. Since $\mathfrak{B}$ is supersymmetric and associative, then $\mathfrak{B}(d,x\mydot y)=\mathfrak{B}(d\mydot x,y)=(-1)^{\vert x\vert + \vert y\vert}\mathfrak{B}(x,y\mydot d)$, which implies that $\Phi(x,y)=-\mathfrak{B}(\delta(x),y)=-(-1)^{\vert x\vert + \vert y\vert}\mathfrak{B}(x,D(y))$. Furthermore,
$\Psi(x)=\mathfrak{B}(x\mydot d,d)=\mathfrak{B}(d,d\mydot x)=-\Theta(x)=\mathfrak{B}(x,\mathfrak{a}_{0})$. Since $\mathcal{A}$ is a left Zinbiel superalgebra satisfying:
$$\mathscr{X}\mydot(\mathscr{Y}\mydot \mathscr{Z})+(-1)^{\vert z\vert (\vert x\vert + \vert y\vert)}(\mathscr{Z}\mydot \mathscr{X})\mydot \mathscr{Y}=0,~~
\mbox{for~all}~\mathscr{X}\in \mathcal{A}_{{\vert x\vert}},
\mathscr{Y}\in \mathcal{A}_{{\vert y\vert}},\mathscr{Z} \in \mathcal{A}_{{\vert z\vert}},$$
then an easy calculation proves that $(\delta,D,\mathfrak{a}_{0})$ is an admissible
triple of $\mathcal{H},$
\begin{center}$\delta(\mathcal{H})\subseteq {\rm Ann}(\mathcal{H}),~\delta(\mathcal{H}^{2})=\{0\},~
\delta\circ D=D\circ \delta=\{0\},~\mathfrak{a}_{0}\in {\rm Ann}(\mathcal{H})\cap \mathcal{H}_{\bar{0}}$ and $\mathfrak{B}(\mathfrak{a}_{0},\mathfrak{a}_{0})=0$.
\end{center} Consequently, we can consider the quadratic Zinbiel superalgebra $\overline{\mathcal{A}}=\mathcal{N}^{*}\oplus \mathcal{H}\oplus \mathcal{N}$, odd double extension of $(\mathcal{H},\cdot,\mathfrak{B}_{\mathcal{H}})$ by $\mathcal{N}$ by means of $(\delta,D,\mathfrak{a}_{0})$.

Clearly, the map
\begin{center}$\Gamma:\mathcal{A}\rightarrow \overline{\mathcal{A}},~\lambda e+x+\lambda'd\longmapsto \lambda d^{*}+x+\lambda'd$,
\end{center} is an isomorphism of Zinbiel superalgebras. Then, $\mathcal{A}$ is isomorphic to the odd double extension of $\mathcal{H}$ by $\mathcal{N}$. \hfill $\square$\\


\end{document}